\documentclass[letterpaper]{amsart} 
 
\usepackage{amssymb,latexsym} 

\numberwithin{equation}{section}

\newtheorem{theorem}{Theorem}[section] 
\newtheorem{proposition}[theorem]{Proposition}

\newtheorem{lemma}[theorem]{Lemma} 
 
\theoremstyle{definition}

\def\ZZ{\mathbb{Z}}

\def\ZZ{\mathbb{Z}}

\begin{document}

\title{Quivers of Finite Mutation Type and Skew-Symmetric Matrices}

\author{Ahmet I. Seven}

\address{Middle East Technical University, 06531, Ankara, Turkey}
\email{aseven@metu.edu.tr}

\thanks{The author's research was supported in part
by Turkish Research Council (TUBITAK)}

\date{April 21, 2010}

\begin{abstract}
Quivers of finite mutation type are certain directed graphs that first arised in Fomin-Zelevinsky's theory of cluster algebras. It has been observed that these quivers are also closely related with different areas of mathematics. In fact, main examples of finite mutation type quivers are the quivers associated with triangulations of surfaces. In this paper, we study structural properties of finite mutation type quivers in relation with the corresponding skew-symmetric matrices. We obtain a characterization of finite mutation type quivers that are associated with triangulations of surfaces and give a new numerical invariant for their mutation classes.
\end{abstract}



\maketitle

\section{Introduction}

\label{sec:introduction}

Quivers of finite mutation type are certain directed graphs that first arised in Fomin-Zelevinsky's theory of cluster algebras. It has been observed that these quivers are also closely related with different areas of mathematics. In fact, main examples of finite mutation type quivers are the quivers associated with triangulations of surfaces as introduced in \cite{FST}. They also provide interesting classes of non-commutative algebras \cite{AB}. A classification of finite mutation type quivers has been obtained recently in \cite{FSTu}. In this paper, we study structural properties of finite mutation type quivers in relation with the corresponding skew-symmetric matrices. We determine a class of subquivers, which we call basic quivers, and show that they have a natural linear-algebraic interpretation. In particular, we obtain a characterization of finite mutation type quivers that are associated with triangulations of surfaces and give a new numerical invariant for their mutation classes. We also give a theoretical proof of the classification of finite mutation type quivers that are not associated with triangulations of a surface (Lemma~\ref{lem:E6-classify}), which was obtained in \cite{FSTu} partly using a computer program. 


To state our results, we need some terminology. Formally, a quiver is a pair $Q=(Q_0,Q_1)$ where $Q_0$ is a finite set of vertices and $Q_1$ is a set of arrows between them. It is represented as a directed graph with the set of vertices $Q_0$ and a directed edge for each arrow. In this paper, we are more  concerned with the number of arrows between the vertices rather than the arrows themselves, so by a quiver we mean a directed graph $Q$, with no loops or 2-cycles, whose edges are weighted with positive integers. If the weight of an edge is $1$, we do not specify it in the picture and call it a single edge; if an edge has weight $2$ we call it a double edge for convenience. If all edges of $Q$ are single edges, we call $Q$ "simply-laced". 
By a {subquiver} of $Q$, we always mean a quiver obtained from $Q$ by taking an induced (full) directed subgraph on a subset of vertices and keeping all its edge weights the same as in $Q$.

For a quiver $Q$ with vertices $1,...,n$, there is the uniquely associated skew-symmetric matrix $B=B^Q$ defined as follows: for each edge $\{i,j\}$ directed from $i$ to $j$, the entry $B_{i,j}$ is the corresponding weight; if $i$ and $j$ are not connected to each other by an edge then $B_{i,j}=0$. Recall from \cite{CAII} that, for each vertex $k$, the mutation of the quiver $Q$ at a vertex $k$ transforms $Q$ to the quiver $Q'=\mu_k(Q)$ whose corresponding skew-symmetric matrix $B'=B^{Q'}$ is the following: $B'_{i,j}=-B_{i,j}$ if $i=k$ or $j=k$; otherwise $B'_{i,j}=B_{i,j}+sgn(B_{i,k})[B_{i,k}B_{k,j}]_+$ (where we use the notation $[x]_+=max\{x,0\}$ and $sgn(x)=x/|x|$ with $sgn(0)=0$). The operation $\mu_k$ is involutive, so it defines a mutation-equivalence relation on quivers (or equivalently on skew-symmetric matrices). A quiver $Q$ is said to be of "finite mutation type" if its mutation-equivalence class is finite. It is well known that, in a finite mutation type quiver with at least three vertices, any edge is a single edge or a double edge; any subquiver is also of finite mutation type. The most basic examples of finite mutation type quivers are Dynkin quivers (Figure~\ref{fig:dynkin-diagrams}), which correspond to skew-symmetric cluster algebras of finite type \cite{CAII}. 

Another important class of finite mutation type quivers has been obtained in \cite{FST} using a construction that associates quivers to certain triangulations of surfaces. 
In this paper, we will not use this construction, so we do not recall it here (we will only use some of their well-known properties).
We call these quivers \emph{quivers that come from the triangulation of a surface}. More recently, it has been shown that these are almost all of the finite mutation type quivers:

\begin{theorem} 
\label{th:fin-mut}\cite[Theorem~6.1]{FSTu} 
A connected quiver $Q $ with at least three vertices is of \emph{finite mutation type} if and only if it comes from the triangulation of a surface or it is mutation-equivalent to one of the exceptional types $E_6$, $E_7$, $E_8$, $E_6^{(1)}$, $E_7^{(1)}$, $E_8^{(1)}$, $E_7^{(1,1)}$, $E_8^{(1,1)}$, $X_6$, $X_7$
(Figures~\ref{fig:dynkin-diagrams},~\ref{fig:exceptional}).
\end{theorem}

\noindent 
The main tool in proving this classification theorem is a purely combinatorial characterization of quivers that come from triangulations of surfaces as quivers that can be composed by matching quivers from a small set of simple quivers. We will not use this construction either, so we do not recall it here. The proof is obtained by determining minimal quivers that are indecomposable, i.e. can not be composed from those simple quivers \cite[Theorem~5.11]{FSTu}.

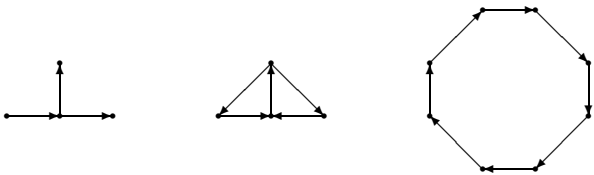
\begin{figure}[ht] 
\[ 
\begin{array}{ccl} 
&& 
\setlength{\unitlength}{1.0pt} 
\begin{picture}(300,60)(0,-2) 
\put(0,0){\circle*{2.0}}
\put(20,0){\circle*{2.0}}
\put(40,0){\circle*{2.0}}
\put(20,20){\circle*{2.0}}

\put(0,0){\vector(1,0){20}}
\put(20,0){\vector(1,0){20}}
\put(20,0){\vector(0,1){20}}

\put(80,0){\circle*{2.0}}
\put(100,0){\circle*{2.0}}
\put(120,0){\circle*{2.0}}
\put(100,20){\circle*{2.0}}

\put(80,0){\vector(1,0){20}}
\put(100,0){\vector(0,1){20}}
\put(100,20){\vector(-1,-1){20}}
\put(100,20){\vector(1,-1){20}}
\put(120,0){\vector(-1,0){20}}

\put(220,0){\circle*{2.0}} 
\put(220,20){\circle*{2.0}}
\put(200,40){\circle*{2.0}}
\put(180,40){\circle*{2.0}}
\put(160,20){\circle*{2.0}}
\put(180,-20){\circle*{2.0}}
\put(160,0){\circle*{2.0}}
\put(200,-20){\circle*{2.0}}

\put(220,0){\vector(-1,-1){20}}
\put(220,20){\vector(0,-1){20}}
\put(200,40){\vector(1,-1){20}}
\put(200,-20){\vector(-1,0){20}}
\put(200,40){\line(-1,0){20}}
\put(180,40){\vector(1,0){20}}
\put(160,20){\vector(1,1){20}}
\put(160,0){\vector(0,1){20}}
\put(180,-20){\vector(-1,1){20}}
\put(180,-20){\line(1,0){20}}

\end{picture} 
\\[.1in] 
\end{array} \] 
\caption{Basic quivers (the cycle has at least four vertices)} 
\label{fig:basic} 
\end{figure}

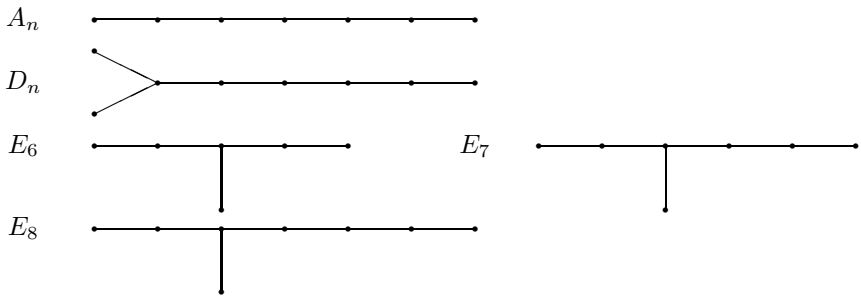
\begin{figure}[ht] 
\vspace{-.2in} 
\[ 
\begin{array}{ccl} 
A_n && 
\setlength{\unitlength}{1.2pt} 
\begin{picture}(300,17)(0,-2) 
\put(0,0){\line(1,0){120}} 
\multiput(0,0)(20,0){7}{\circle*{2}} 



\end{picture}\\
D_n 
&& 
\setlength{\unitlength}{1.2pt} 
\begin{picture}(140,17)(0,-2) 
\put(20,0){\line(1,0){100}} 
\put(0,10){\line(2,-1){20}} 
\put(0,-10){\line(2,1){20}} 
\multiput(20,0)(20,0){6}{\circle*{2}} 
\put(0,10){\circle*{2}} 
\put(0,-10){\circle*{2}} 
\end{picture} 
\\
E_6 
&& 
\setlength{\unitlength}{1.2pt} 
\begin{picture}(240,17)(0,-2) 
\put(0,0){\line(1,0){80}} 
\put(40,0){\line(0,-1){20}} 
\put(40,-20){\circle*{2}} 
\multiput(0,0)(20,0){5}{\circle*{2}}

\put(140,0){\line(1,0){100}} 
\put(180,0){\line(0,-1){20}} 
\put(180,-20){\circle*{2}} 

\put(140,0){\circle*{2}} 
\put(160,0){\circle*{2}} 
\put(180,0){\circle*{2}} 
\put(200,0){\circle*{2}} 
\put(220,0){\circle*{2}} 
\put(240,0){\circle*{2}} 
\put(120,0){\makebox(0,0){$E_7$}}

\end{picture} 
\\[.1in] 
E_8 
&& 
\setlength{\unitlength}{1.2pt} 
\begin{picture}(140,17)(0,-2) 
\put(0,0){\line(1,0){120}} 
\put(40,0){\line(0,-1){20}} 
\put(40,-20){\circle*{2}} 
\multiput(0,0)(20,0){7}{\circle*{2}} 
\end{picture} 
\\[.2in] 
\end{array} 
\] 
\vspace{-.1in} 
\caption{Dynkin quivers: each edge is assumed to be arbitrarily oriented} 
\label{fig:dynkin-diagrams} 
\end{figure}

\begin{figure}[ht] 
\vspace{-.2in} 
\[ 
\begin{array}{ccl} 
E_6^{(1)} 
&& 
\setlength{\unitlength}{1.2pt} 
\begin{picture}(260,37)(0,-2) 
\put(0,0){\line(1,0){80}} 
\put(40,0){\line(0,-1){40}} 
\put(40,-20){\circle*{2}}
\put(40,-40){\circle*{2}}
\multiput(0,0)(20,0){5}{\circle*{2}}

\put(140,0){\line(1,0){120}} 
\put(200,0){\line(0,-1){20}} 
\put(200,-20){\circle*{2}} 
\put(260,0){\circle*{2}}

\put(140,0){\circle*{2}} 
\put(160,0){\circle*{2}} 
\put(180,0){\circle*{2}} 
\put(200,0){\circle*{2}} 
\put(220,0){\circle*{2}} 
\put(240,0){\circle*{2}} 
\put(120,0){\makebox(0,0){$E_7^{(1)}$}}

\end{picture} 
\\[.5in] 
E_8^{(1)} 
&& 
\setlength{\unitlength}{1.2pt} 
\begin{picture}(160,17)(0,-2) 
\put(0,0){\line(1,0){140}} 
\put(40,0){\line(0,-1){20}} 
\put(40,-20){\circle*{2}} 
\put(140,0){\circle*{2}} 
\multiput(0,0)(20,0){7}{\circle*{2}} 
\end{picture} 
\\[.3in] 

E_6^{(1,1)}
&& 
\setlength{\unitlength}{1.2pt} 
\begin{picture}(140,37)(0,-2) 

\put(0,0){\line(1,0){20}} 
\put(120,0){\line(1,0){20}} 
\put(60,0){\line(1,0){20}} 

\put(40,-20){\vector(0,1){40}}
\put(40,20){\vector(1,-1){20}}
\put(20,0){\vector(1,-1){20}}

\put(40,20){\vector(-1,-1){20}}
\put(60,0){\vector(-1,-1){20}}

\put(40,20){\vector(4,-1){80}}
\put(120,0){\vector(-4,-1){80}}

\put(40,-20){\circle*{2}} 
\put(40,20){\circle*{2}}

\multiput(0,0)(20,0){2}{\circle*{2}} 
\multiput(60,0)(20,0){2}{\circle*{2}} 
\multiput(120,0)(20,0){2}{\circle*{2}} 

\put(43,0){\makebox(0,0){$2$}}

\end{picture} 
\\[.3in] 

E_7^{(1,1)}
&& 
\setlength{\unitlength}{1.2pt} 
\begin{picture}(140,37)(0,-2) 

\put(0,0){\line(1,0){40}} 
\put(140,0){\line(1,0){40}} 

\put(60,-20){\vector(0,1){40}}
\put(60,20){\vector(1,-1){20}}
\put(40,0){\vector(1,-1){20}}

\put(60,20){\vector(-1,-1){20}}
\put(80,0){\vector(-1,-1){20}}

\put(60,20){\vector(4,-1){80}}
\put(140,0){\vector(-4,-1){80}}

\put(60,-20){\circle*{2}} 
\put(60,20){\circle*{2}}

\multiput(0,0)(20,0){3}{\circle*{2}} 
\multiput(80,0)(20,0){1}{\circle*{2}} 
\multiput(140,0)(20,0){3}{\circle*{2}} 

\put(63,0){\makebox(0,0){$2$}}

\end{picture} 
\\[.2in] 

E_8^{(1,1)}
&& 
\setlength{\unitlength}{1.2pt} 
\begin{picture}(140,37)(0,-2) 

\put(0,0){\line(1,0){20}} 
\put(120,0){\line(1,0){80}} 

\put(40,-20){\vector(0,1){40}}
\put(40,20){\vector(1,-1){20}}
\put(20,0){\vector(1,-1){20}}

\put(40,20){\vector(-1,-1){20}}
\put(60,0){\vector(-1,-1){20}}

\put(40,20){\vector(4,-1){80}}
\put(120,0){\vector(-4,-1){80}}

\put(40,-20){\circle*{2}} 
\put(40,20){\circle*{2}}

\multiput(0,0)(20,0){2}{\circle*{2}} 
\multiput(60,0)(20,0){1}{\circle*{2}} 
\multiput(120,0)(20,0){5}{\circle*{2}} 

\put(43,0){\makebox(0,0){$2$}}

\end{picture} 
\\[.2in] 
X_6 
&& 
\setlength{\unitlength}{1.2pt} 
\begin{picture}(260,42)(0,-2) 
\put(20,0){\vector(-1,0){20}} 
\put(0,0){\vector(0,1){20}} 
\put(0,20){\vector(1,-1){20}} 
\multiput(20,0)(20,0){2}{\circle*{2}} 
\put(0,20){\circle*{2}} 
\put(0,0){\circle*{2}} 
\put(40,20){\circle*{2}} 
\put(20,-20){\circle*{2}} 

\put(20,0){\line(0,-1){20}} 
\put(20,0){\vector(1,0){20}} 
\put(40,0){\vector(0,1){20}} 
\put(40,20){\vector(-1,-1){20}} 

\put(-3,7){\makebox(0,0){$2$}} 
\put(43,7){\makebox(0,0){$2$}}

\put(180,0){\vector(-1,0){20}} 
\put(160,0){\vector(0,1){20}} 
\put(160,20){\vector(1,-1){20}} 

\put(160,20){\circle*{2}} 
\put(160,0){\circle*{2}} 
\put(200,20){\circle*{2}} 
\put(180,-20){\circle*{2}} 
\put(180,0){\circle*{2}} 
\put(200,0){\circle*{2}} 

\put(180,0){\vector(0,-1){20}} 
\put(180,-20){\vector(1,0){20}} 
\put(200,-20){\vector(-1,1){20}} 

\put(180,0){\vector(1,0){20}} 
\put(200,0){\vector(0,1){20}} 
\put(200,20){\vector(-1,-1){20}} 

\put(157,7){\makebox(0,0){$2$}} 
\put(203,7){\makebox(0,0){$2$}} 

\put(130,2){\makebox(0,0){$X_7$}}
\put(190,-24){\makebox(0,0){$2$}}

\end{picture} 
\\[.4in] 
\end{array} 
\] 
\vspace{-.1in} 
\caption{exceptional quivers of finite mutation type which are not Dynkin: edges with unspecified orientation are assumed to be arbitrarily oriented} 
\label{fig:exceptional} 
\end{figure}
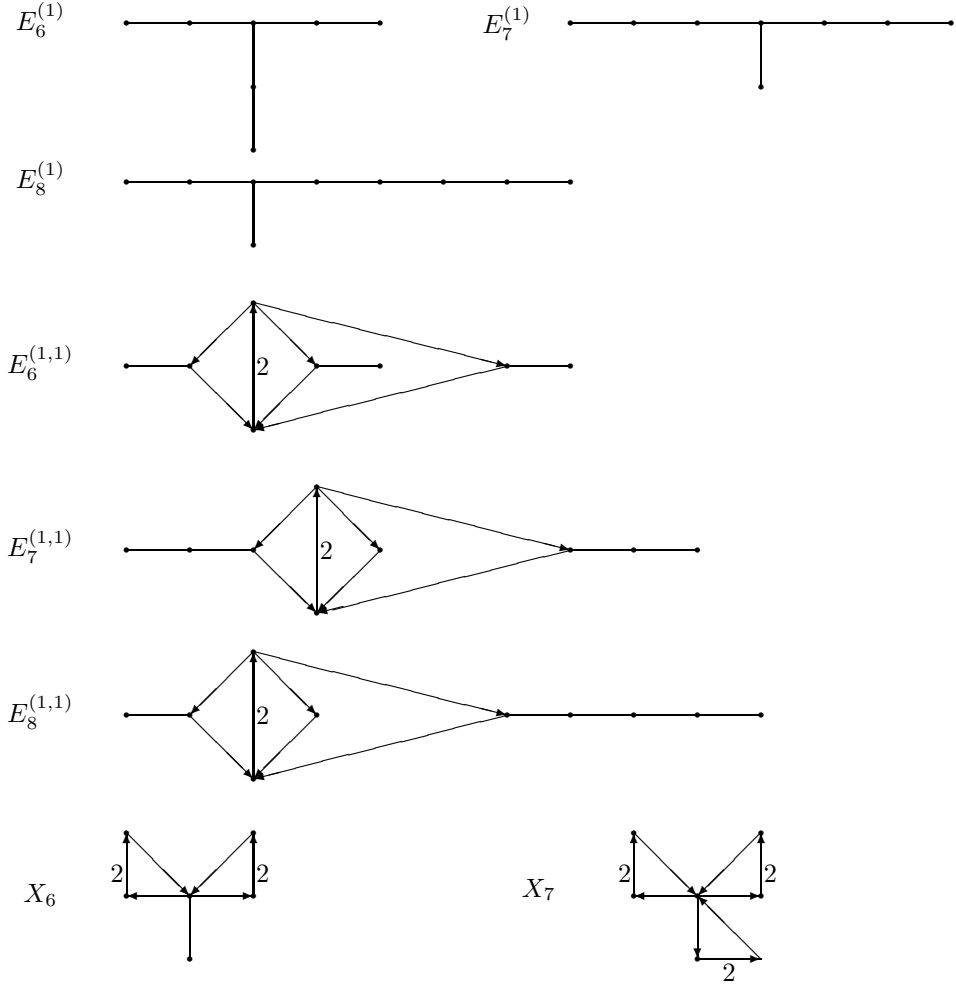

In this paper, to understand the structure of finite mutation type quivers, we identify another class of subquivers that we call "basic (sub)quivers" and use them give an algebraic/combinatorial characterization of the finite mutation type quivers that come from triangulations of surfaces. More explicitly, we define a basic quiver as one of the following: a Dynkin tree $D_4$, two adjacent oriented simply-laced triangles, an oriented simply-laced cycle with at least four vertices (see Figure~\ref{fig:basic}). Here by a cycle we mean a subquiver whose vertices can be labeled by elements of $\ZZ/m\ZZ$ so that the edges betweeen them are precisely $\{i,i+1\}$ for $i \in  \ZZ/m\ZZ$. To proceed, we need a little bit more terminology. For each vertex $i$ in a quiver $Q$ with vertex set $\{1,2,...,n\}$, we denote by $e_i$ the $i$-th standard basis vector of $\ZZ^n$. For any vector $u$ in $\ZZ^n$, we define $supp_Q(u)$ to be the subquiver of $Q$ on the vertices which correspond to the non-zero coordinates of $u$ and call it the support of $u$ in $Q$. Now we can state our first main result:

\begin{theorem}\label{th:basic radical}
Suppose that $Q$ is a finite mutation type quiver with at least three vertices. Then $Q$ comes from the triangulation of a surface
if and only if the following holds for any basic subquiver $S$:
\begin{enumerate}
\item[(i)] 
$S$ contains a subquiver of the form $supp_Q(u)$ where $u$ is a non-zero radical vector of $B^Q$ such that each non-zero coordinate of $u$ is either $1$ or $-1$ (here $u$ is radical if $B^Qu=0$). The subquiver $supp_Q(u)$ has exactly two vertices or it is a cycle.
\item[(ii)] 
furthermore if $S$ is an oriented cycle of length at least $5$, then the vector $u$ whose coordinates corresponding to the vertices of $S$ is $1$ and $0$ in the remaining vertices is a radical vector for $B_Q$ (in particular $S=supp_Q(u)$).

\end{enumerate}
\end{theorem}

Thus we have, in particular, obtained an algebraic interpretation of basic subquivers in quivers that come from the triangulation of a surface. 
We will also give a numerical invariant for their mutation classes which is related to this interpretation, involving another common class of subquivers as well: double edges and non-oriented cycles. For this purpose, it turns out to be convenient to work in $\bar{V}:=\ZZ^n/2\ZZ^n$, which is a vector space over $\ZZ/2\ZZ$ (which is the field with two elements). To be more precise, for a finite mutation type quiver $Q$, 
we denote by $\bar{B}^Q$ the skew-symmetric matrix whose entries are the corresponding entries of ${B}^Q$ modulo $2\ZZ$. We denote by $\bar{V}^{Q}_0$ the space of radical vectors of $\bar{B}^Q$ (over $\ZZ/2\ZZ$); we call a vector $u$ in $\bar{V}^{Q}_0$ a "basic radical vector" if $supp_Q(u)$ has exactly two vertices or it is a cycle (oriented or not). We denote by $\bar{V}^{Q}_{00}$ the subspace spanned by the basic radical vectors of $\bar{B}^Q$ over $\ZZ/2\ZZ$; if there are no basic radical vectors, then we take $\bar{V}^{Q}_{00}$ as the zero subspace. Let us also note that the radical vectors given by Therem~\ref{th:basic radical} are basic radical vectors for $\bar{B}^Q$. Our next result relates these vectors to subquivers: 

\begin{theorem}
\label{th:V00}
Suppose that $Q$ is a connected finite mutation type quiver with at least three vertices. Suppose also that $S$ is a subquiver which is a double edge or a non-oriented cycle. Let $u$ be the vector whose coordinates corresponding to the vertices of $S$ is $1$ and $0$ in the remaining vertices. Then $u$ is a radical vector for $\bar{B}_Q$.

Furthermore, if $Q$ comes from the triangulation of a surface or it is mutation-equivalent to one of $X_6,X_7$, then we have the following: 
\begin{enumerate}
\item[(i)] 
$dim(\bar{V}^{Q}_{0}/\bar{V}^{Q}_{00})\leq 1$. 
\item[(ii)] 
if $Q$ and $Q'$ are mutation-equivalent, then $dim(\bar{V}^{Q}_{00})=dim(\bar{V}^{Q'}_{00})$.
\end{enumerate}
\end{theorem}

\noindent
Let us note, in particular, that $dim(\bar{V}^{Q}_{00})$ is a numerical invariant for the mutation classes of quivers that come from triangulation of a surface. In view of Theorem~\ref{th:basic radical}, it can be considered as a count of subquivers $S$ such that $S$ is a double edge or a non-oriented cycle or a basic quiver, modulo those which overlap in a way that the supports of the corresponding basic radical vectors coincide.
Let us also note that the first part of the theorem holds for any finite mutation type quiver. However, the second part may not be true for a quiver which belongs to one of the types $E$ in Theorem~\ref{th:fin-mut}. Also, in part (i), the equality may hold; e.g. it holds for a Dynkin quiver $Q$ which is of type $A_{2n+1}$, $n\geq 2$. 




In the classification theorem \cite[Theorem~6.1]{FSTu}, which is Theorem~\ref{th:fin-mut} above, the classification of quivers that do \emph{not} come from the triangulation of a surface was done, in part, using a computer program \cite[Proof of Theorem 6.1]{FSTu}. Here, using our approach, we suggest an algebraic/combinatorial proof (Lemma~\ref{lem:E6-classify}). More precisely, we show the following: 

\begin{theorem}
\label{th:E6 inv}
Let $Q$ be a connected quiver of finite mutation type. Suppose also that $Q$ has a subquiver which is mutation-equivalent to $E_6$ (resp. $X_6$). Then any quiver which is mutation-equivalent to $Q$ also contains a subquiver which is mutation-equivalent to $E_6$ (resp. $X_6$). Furthermore $Q$ is mutation-equivalent to a quiver which is one of the types $E$ (resp. $X$) given in Theorem~\ref{th:fin-mut}. 
\end{theorem}




We prove our results in Section~\ref{sec:proof} after some preparation in Section~\ref{sec:pre}.













\section{Preliminary results}
\label{sec:pre}

In this section, we will recall some more terminology and prove some statements that we will use to prove our results. As we discussed in Section~\ref{sec:introduction}, we mainly study quivers of finite mutation type. There is also a stronger notion of \emph{finite (cluster) type}: a quiver $Q$ is said to be of finite type if any edge in any quiver $Q'$ which is mutation-equivalent to $Q$ is a single edge. Quivers of finite type were classified by Fomin and Zelevinsky in \cite{CAII}. Their classification is identical to the Cartan-Killing classification. More precisely, a quiver $Q$ is of finite type if and only if $Q$ is mutation-equivalent to an orientation of a Dynkin quiver (Figure~\ref{fig:dynkin-diagrams}). Another related definition is the following: a quiver $Q$ is said to be of \emph{minimal infinite} type if it is of infinite type and any proper subquiver of $Q$ is of finite type. A list of minimal infinite type quivers has been obtained in \cite{S2}. In particular, any minimal infinite type quiver with at least three vertices is mutation-equivalent to an extended Dynkin quiver \cite[Theorem~3.2]{S2}.

Let us also recall that the mutation operation can be viewed as a base change transformation for a skew-symmetric bilinear form \cite[Section~2]{FSTu}. To be more specific, let $Q$ be a quiver with vertices $1,...,n$. Let $\Omega$ be the bilinear form on $\ZZ^n$ defined as follows: $\Omega(e_i,e_j)=B^{Q}_{i,j}$; here $e_i$ denotes the $i$-th standard basis vector. Then $\mu_k(B^{Q})$ is the matrix that represents $\Omega$ with respect to the basis $\{e'_j\}$: $e'_k=-e_k$, $e'_j=e_j$ if $\Omega(e_k,e_j)>0$, $e'_j=e_j-\Omega(e_k,e_j)e_k$ if $\Omega(e_k,e_j)<0$. The coordinates of a vector $u=(u_1,...,u_n)$ in $\ZZ^n$ with respect to this new basis will be the same except in the $k-$th coordinate, which becomes $-u_k+\sum u_j$ over all $j$ such that $\Omega(e_k,e_j)<0$. In particular, the rank of $B^{Q}$ is invariant under the mutation operation.

Finally, let us recall that a vertex $i$ in a quiver is called a {source} ({rsp. sink}) if all adjacent edges are oriented away (resp. towards) $i$. 
A quiver is called {acyclic} if it has no oriented cycles at all. It is well known that an acyclic quiver has a source and a sink. Also acyclic quivers  of finite mutation type with at least three vertices are the Dynkin and extended Dynkin quivers \cite{BR}.

Let us now give some basic examples of quivers which are of infinite mutation type (i.e. not of finite mutation type):
\begin{proposition}\label{prop:infinite}
Let $Q$ be a connected quiver which has at least three vertices.
\begin{enumerate}
\item[(i)] If $Q$ is of finite mutation type, then any edge of $Q$ is a single edge or a double edge. 

\item[(ii)] Suppose that $Q$ has exactly three vertices and has a double edge. Then $Q$ is of finite mutation type if and only if it is an oriented triangle with edge weights $2,1,1$ or $2,2,2$. 

\item[(iii)] Any non-simply-laced, non-oriented cycle is of infinite mutation type.

\item[(iv)] Suppose that $Q$ is a simply-laced quiver.  
If $Q$ contains a non-oriented cycle $C$ such that there is a vertex $k$ which is connected to exactly an odd number of vertices in $C$, then $Q$ is of infinite mutation type. If $C$ is an oriented cycle in $Q$ and $k$ is connected to exactly an odd number greater than or equal to $3$ vertices in $C$, then $Q$ is also of infinite mutation type.

\item[(v)] Suppose that $Q$ is a quiver which has no oriented cycles but has at least two non-oriented cycles. Then $Q$ is of infinite mutation type. 
\end{enumerate}
\end{proposition}

\noindent
\begin{proof}
Statements (i),(ii),(iii) are obtained easily from the definitions by observing that if the conclusions do not hold then, by an iterative process, the mutation class of the quiver contains edges of arbitrarily large weights. Let us now prove (iv) for a non-oriented cycle $C$. (The second part for an oriented cycle  follows by similar arguments). By part (iii), we can assume that $C$ is simply-laced. First we consider the case where $k$ is connected to exactly one vertex, say $c$, in $C$. Let us suppose first that $C$ is a triangle. Applying a mutation at a source or sink of $C$ if necessary, we can assume that $c$ is a source or sink; mutating at the vertex which is neither a source or sink, we obtain a quiver which contains a three-vertex tree which has a double edge; then part (ii) applies. Let us now suppose that $C$ has more than $3$ vertices. Then, applying a mutation at a source or sink of $C$ if necessary, we can assume that there is a vertex $c'\ne c$ in $C$ which is neither a source nor a sink in $C$. Then, in $\mu_{c'}(Q)$, the subquiver $C'$ obtained from $C$ by removing $c'$ is a non-oriented cycle and $k$ is connected to exactly one vertex in $C'$. Then the statement (iv) follows by induction.

Let us now consider the case where $k$ is connected to exactly three vertices in $C$. Then there are three cycles, say $C_1,C_2,C_3$, that contain $k$; one of them, say $C_1$, is necessarily non-oriented. 
If one of the cycles $C_2$ or $C_3$ has more than three vertices, then there is a vertex in that cycle connected to exactly one vertex in $C_1$, which is the case we considered above. Thus we can further assume that $C_2$ and $C_3$ are triangles. 
Given all this, we proceed as follows. If $C$ has exactly three vertices, then the statement follows from a direct check. 
If $C$ has more than three vertices, then one of the cycles $C_1,C_2,C_3$ also has more than three vertices; since $C_2$ and $C_3$ are triangles, the cycle $C_1$ must have at least four vertices. If any of $C_2$ or $C_3$ is non-oriented, then there is a vertex in $C_1$ which is connected to exactly one vertex in that cycle, which is the case we considered above. Then the only subcase left to consider is the case where both $C_2$ and $C_3$ are oriented. Then, in $\mu_{k}(Q)$, the subquiver $C\cup \{k\}$ consists of a non-oriented cycle $C'$ that contains $k$ and an additional vertex which is connected to exactly one vertex in $C'$, which is again the case we considered above. To consider the case where $k$ is connected to at least five vertices in $C$, we note that in this case there is a non-oriented cycle $C'$ which contains $k$ and there is a vertex in $C$ connected to exactly one vertex in $C'$, which is also the case we considered.


To prove part (v), we can assume that any cycle in $Q$ is simply-laced by part(iii). Let us now suppose that $C$ is a cycle in $Q$ with minimal number of vertices. There is a vertex $k$ which is not in $C$ but connected to $C$. If $k$ is connected to $C$ by a double edge $e$, then there is a three-vertex acyclic subquiver that contains $e$, so part (ii) applies. Thus we can also assume that any edge connecting $k$ to $C$ is a single edge. If $k$ is connected to an odd number of vertices in $C$, then part (iv) applies. If $k$ is connected to an even number of vertices and $C$ is a triangle or a square, then the statement follows from a direct check; if $C$ has at least five vertices, then there is a non-oriented cycle $C'$ containing $k$ such that there is a vertex $r\ne k$ which is connected to exactly an odd number of vertices in $C'$, so part (iv) applies. This completes the proof of the proposition.
\end{proof}


The following statement follows from the definitions:
\begin{lemma}\label{lem:Srad}
Suppose that $Q$ is a quiver and $S=\{s_1,...,s_r\}$ is a subquiver of $Q$. Let $u=e_{s_1}+...+e_{s_r}$. Then we have the following:
\begin{enumerate}
\item[(i)]
$u$ is a radical vector for $B^Q$ if and only if, for any vertex $j$ and for the edges connecting $j$ to (a vertex in) $S$, the following holds: the number of such edges entering $j$ is equal to the number of the ones leaving, each edge being multipled by its weight. 
\item[(ii)]
$u$ is a radical vector for $\bar{B}^Q$ if and only if, for any vertex $j$, the sum of the weights of the edges connecting $j$ to $S$ is even.
\end{enumerate}
\end{lemma}

We also give several different characterizations of the mutation class of the Dynkin quiver $E_6$:
\begin{proposition}\label{prop:E6 char}
Suppose that $Q$ is a simply-laced connected quiver with six vertices such that $Q$ does not contain any non-oriented cycle. Then the following statements (a),(b),(c),(d) are equivalent.

\begin{enumerate}
\item[(a)] $Q$ is mutation-equivalent to $E_6$.
\item[(b)] $Q$ contains a basic subquiver and $B^Q$ has corank $0$.
\item[(c)] $Q$ contains a basic subquiver and $\bar{B}^Q$ has corank $0$.
\item[(d)] The following (i),(ii),(iii) hold:

\begin{enumerate}
\item[(i)] $Q$ contains a basic subquiver,
\item[(ii)] for each cycle $C$ in $Q$, there is a vertex which is connected to exactly 
one vertex in $C$, 
\item[(iii)] for each pair of vertices which are not connected to each other, there is a vertex which is connected to exactly one of them.

\end{enumerate}

\end{enumerate}


\end{proposition}



\noindent
We will first prove that (a) and (b) are equivalent. To show that (a) implies (b), let us suppose that $Q$ is mutation-equivalent to $E_6$. 
If $Q$ does not contain any basic quiver, then it is mutation-equivalent to the Dynkin quiver $A_6$ \cite[Corollary~5.15]{S2}. Therefore $Q$ contains a basic subquiver. It also follows from a direct computation that $B^Q$ has corank $0$ (it is enough to compute it for $Q=E_6$ because the rank of $B^Q$ is invariant under the mutation operation). For the converse, let us suppose that $B^Q$ has corank $0$ and $Q$ contains a basic subquiver $S$. Let us first assume that $Q$ is of finite (cluster) type. Then, by the classification of finite type quivers, $Q$ is mutation-equivalent to a Dynkin quiver which is of type $A_6$ or $D_6$ or $E_6$. Since $Q$ contains a basic subquiver, it is not mutation-equivalent to $A_6$; since $B^Q$ has corank $0$, the quiver $Q$ is not mutation-equivalent to $D_6$ either. Therefore $Q$ is mutation-equivalent to $E_6$. Let us now assume that $Q$ is not of finite type. Then $Q$ contains a minimal infinite type subquiver $M$. Then $M$ is mutation-equivalent to the extended Dynkin tree $D^{(1)}_4$ \cite[Theorem~3.2]{S2}. The skew-symmetric matrix $B^M$ has  corank $3$, so the corank of $B^Q$ is at least $2$, contradicting our assumption. Thus $Q$ is necessarily mutation-equivalent to $E_6$. This completes the proof of the equivalence of (a) and (b). Using similar arguments, the equivalence of (a) and (c) can be proved easily. The equivalence of the other statements follow by Lemma~\ref{lem:Srad} and Proposition~\ref{prop:infinite}(iv).

 
Let us also record the following statement which follows immediately from the previous proposition and Proposition~\ref{prop:infinite}(iv): 
\begin{proposition}\label{prop:E6 reorient}
Suppose that $Q$ is a simply-laced quiver which contains a non-oriented cycle. If the underlying (undirected) graph satisfies (i),(ii,(iii) of Proposition~\ref{prop:E6 char}, then it is of infinite mutation type.

In particular, if the underlying (undirected) graph of $Q$ is equal to the underlying graph of a quiver which is mutation-equivalent to $E_6$, then $Q$ is of infinite mutation type.
\end{proposition}




\section{Proofs of Main Results}
\label{sec:proof}

\subsection{Proof of Theorem~\ref{th:basic radical}}

\label{sec:basic radical}

We first show the "only if part", i.e. if $Q$ comes from the triangulation of a surface then it satisfies (i),(ii). For this it is enough to establish the theorem for $Q$ which does not contain any subquiver which is mutation-equivalent to $E_6$ or $X_6$ (because quivers that come from the triangulation of a surface have this property \cite[Corollary~5.13]{FSTu}). We show this by induction on the number of vertices of $S$. The basis of the induction is for $S$ with exactly four vertices. There are three types of such basic quivers: Dynkin tree $D_4$, two adjacent oriented triangles or oriented square (all are simply-laced). For convenience, we will first prove for $S$ which is an oriented square. 

Let us now assume that $S=\{s_1,s_2,s_3,s_4\}$ is oriented cyclically (so $s_1\to s_2 \to s_3 \to s_4 \to s_1$). We will show that one of the vectors $e_{s_1}+e_{s_2}+e_{s_3}+e_{s_4}$ or $e_{s_1}+e_{s_3}$ or $e_{s_2}+e_{s_4}$ is a radical vector using Lemma~\ref{lem:Srad}(i). This trivially follows from Lemma~\ref{lem:Srad} if there is no vertex which is connected to $S$. Thus we can assume that there is a vertex $k$, which is not in $S$, connected to (at least one vertex in) $S$. For any such $k$, we denote the subquiver $\{S,k\}$ by $Sk$ for convenience. If $k$ is connected to $S$ by a double edge, then there is necessarily a three-vertex subquiver which is not as in Proposition~\ref{prop:infinite}(ii), so we assume that \emph{any edge connecting a vertex $k$ to $S$ is a single edge}. Below we will establish the theorem considering possible cases for $k$ to connect to $S$. During the analysis, if we do not specify an orientation on a subquiver, we assume that it is oriented as required by Proposition~\ref{prop:infinite} to be of finite mutation type.

We first show the theorem in the case that (*) for any $k$ which is connected to $S$, the quiver $Sk$ does not contain any non-oriented cycle. Then there are the following two subcases: (i) $k$ is connected to exactly one vertex in $S$ or (ii) $k$ is connected to exactly two vertices $s_i,s_j$ and  the subquiver $\{k,s_i,s_j\}$ is an oriented triangle. 
For the subcase (i) let us assume without loss of generality that $k$ is connected to $s_1$.  We show that $e_{s_2}+e_{s_4}$ is a radical vector. This is true if no vertex (outside $S$) is connected to $s_2$ or $s_4$. Similarly, it is also true if any vertex which is connected to $s_2$ is connected to $s_4$ with the opposite orientation (recall that any edge incident to any $s_i$ is a single-edge). We now consider the remaining two possibilities: (a) There is a vertex $r$ which is connected to exactly one of $s_2,s_4$, say connected to $s_2$ (by a single edge). Note that by asumption (*), the subquiver $Sr$ does not contain any non-oriented cycles either (in particular, $r$ is connected to at most one of $s_1,s_3$). If $r$ is connected to $k$ by a double edge then, the subquiver $\{k,s_2,r\}$ is a three-vertex tree, so contradiction by Proposition~\ref{prop:infinite}(iv). Thus we can assume now that $Skr$ is simply-laced. Then we have the following: If $r$ is not connected to both of $k$ and $s_3$, then the subquiver $Skr$ is mutation-equivalent to $E_6$ (recall our convention that the edges adjacent to $k$ or $r$ are oriented as required by Proposition~\ref{prop:infinite}); if $r$ is connected to both of $k$ and $s_3$, then the subquiver $\{k,r,s_3,s_4,s_1\}$ is a cycle and $s_2$ is connected to exactly three vertices there so contradiction by Proposition~\ref{prop:infinite}(iv). (b) There is a vertex $r$ which is connected to both of $s_2,s_4$ with the same orientation. Then the subquiver $Sr$ contains a non-oriented cycle, contradicting (*).

For the subcase (ii), let us assume without loss of generality that $k$ is connected to $s_1$ and $s_2$ (such that the triangle $\{k,s_1,s_2\}$ is oriented). We show that $e_{s_1}+e_{s_2}+e_{s_3}+e_{s_4}$ is a radical vector. If this is not true then, by Lemma~\ref{lem:Srad}, there is a vertex $r$ as in the following two subcases: (a) $r$ is connected to exactly an odd number of vertices in $S$. Then, by Proposition~\ref{prop:infinite}(iv) or (*), the vertex $r$ is connected to exactly one vertex in $S$. This is the same situation as we considered in the previous subcase (exchanging $r$ and $k$), which implies that $e_{s_1}+e_{s_3}$ or $e_{s_2}+e_{s_4}$ is radical, however this is not true in this case (Lemma~\ref{lem:Srad}). 
(b) $r$ is connected to exactly an even number of vertices in $S$ such that the number of corresponding edges which enter $r$ is different from the ones which leave. Then the subquiver $Sr$ contains a non-oriented cycle, contradicting (*).



Thus for the rest of the proof, we can assume that \emph{there is a vertex $k$ which is contained in a non-oriented cycle $C\subset Sk$}. Then $k$ is connected to at least two vertices in $S$. If $k$ is connected to exactly two vertices $s_i,s_j$ in $S$ and $s_i,s_j$ are connected, then $C$ is the triangle $\{k,s_i,s_j\}$ and one of the remaining vertices in $S$ is connected to exactly one vertex in $C$, which implies that $Q$ is not of finite mutation type (Proposition~\ref{prop:infinite}(iv)), contradicting our assumption. We have the same contradiction if $k$ is connected to exactly three vertices in $S$. Therefore, we can assume that $k$ is connected to exactly two vertices in $S$ which are not connected to each other or $k$ is connected to all four vertices. We proceed considering possible (sub)cases:

{Case 1.} \emph{$k$ is connected to exactly two vertices, say $s_1,s_3$, in $S$, which are not connected to each other.}

{Subcase 1.1.} \emph{$k$ is a source or sink in $C$}. Then $Sk$ has two non-oriented cycles $C=\{k,s_1,s_2,s_3\}$ and $C'=\{k,s_1,s_3,s_4\}$ . We will show that $e_{s_2}+e_{s_4}$ is radical. Suppose that this is not true. Then there is a vertex $r$ which is (i) connected to exactly one of $s_2,s_4$ or is (ii) connected to both of $s_2,s_4$ with the same orientation (Lemma~\ref{lem:Srad}).

In the former case (i), assume without loss of generality that $r$ is connected to $s_2$ (by a single edge) and not connected to $s_4$ (note that $r$ may be connected to $s_1$ or $s_3$). If $r$ is connected to $k$ by a double edge, then the subquiver $\{r,k,s_2\}$ is a three-vertex tree that contains this double edge, contradicting the assumption that $Q$ is of finite mutation type by Proposition~\ref{prop:infinite}(ii). 
If $r$ is connected to $k$ by a single edge, then the vertex $r$ is connected to exactly an odd number of vertices in $C$ or $C'$ (more explicitly, if $r$ is connected to an even number of vertices in $C'$ then it is connected to exactly one of $s_1,s_3$, because $r$ is not connected to $s_4$, then $r$ is connected to exactly three vertices in $C$), contradiction by Proposition~\ref{prop:infinite}(iv). Similarly, if $r$ is not connected to $k$ then the vertex $r$ is connected to exactly an odd number of vertices in $C$ or $C'$ (more explicitly, if $r$ is not connected to any of $s_1,s_3$, then it is connected to exactly one vertex, which is $s_2$, in $C$; otherwise $r$ is connected to exactly one of $s_1,s_3$, then $r$ is connected to exactly one vertex in $C$), contradiction.



 
In the latter case (ii), first we note that $r$ is connected to $k$ because otherwise there is a non-oriented cycle $C''$ that contains $r$  
(because $S$ is oriented) such that $k$ is connected to exactly one vertex in $C''$, contradiction by 
Proposition~\ref{prop:infinite}(iv). If $r$ is connected to $k$ by a double edge, then there is the three vertex tree $\{k,r,s_2\}$, contradiction by Proposition~\ref{prop:infinite}(ii). If $r$ is connected to $k$ by a single edge and not connected to any of $s_1,s_3$, then the subquiver $Skr$ is not of finite mutation type (Proposition~\ref{prop:infinite}(v)) or it is mutation-equivalent to $X_6$; if $r$ is connected to $k$ (by a single edge) and connected to both of $s_1,s_3$, then $Skr$ is mutation-equivalent to $X_6$, contradiction. If $r$ is connected to exactly one of $s_1,s_3$, then it is connected to to exactly three vertices in $S$, so Proposition~\ref{prop:infinite}(iv) applies.


{Subcase 1.2.} \emph{k is a not a source or sink in $C$.} Then $Sk$ consists exactly of an oriented cycle, say $C=\{s_1,s_2,s_3,k\}$, and a non-oriented cycle $C'=\{k,s_1,s_3,s_4\}$ (both containing $k$). 
Note then that any edge incident to $k$ is a single edge by Proposition~\ref{prop:infinite}(ii).

To proceed let us first note the following:

Claim: If a vertex $r$ which is not in $Sk$ is connected to $k$ (by a single edge), then $r$ is connected to exactly one of $s_1,s_3$. 

Proof: if $r$ is not connected to any of them, then it is connected to $s_4$ (by Proposition~\ref{prop:infinite}(iv) because $C'$ is non-oriented), which implies that both cycles  $\{k,s_1,s_4,r\}$, $\{k,r,s_4,s_3\}$ are non-oriented (because there $s_1$ or $s_3$ is a sink or source respectively), then Proposition~\ref{prop:infinite}(v) applies to their union to give a contradiction. If $r$ is connected to both of  $s_1,s_3$, then $r$ is connected to $s_4$ as well (to connect to an even number of vertices in $C'$), so there are four triangles in $C'r$. By Proposition~\ref{prop:infinite}(v), we can assume that exactly two of them, say $T,T'$, are non-oriented and they are not adjacent (it is not possible that all four of these triangles are oriented because $S$ is oriented). 
We may assume, without loss of generality, that $T$ does not contain $k$.
Then we have the following: if $r$ is not connected to $s_2$, then $s_2$ is connected to exactly one vertex in $T$ 
, so Proposition~\ref{prop:infinite}(iv) applies to give a contradiction; if $r$ is connected to $s_2$, then  each of the triangles $T_1=\{r,s_1,s_2\}$ and  $T_2=\{r,s_3,s_2\}$ is non-oriented 
(more explicitly, e.g., if the triangle $T_1$ is oriented, then the triangles $\{r,s_1,k\}$ and $\{r,s_1,s_4\}$ are both non-oriented and adjacent, so Proposition~\ref{prop:infinite}(v) applies), furthermore one of $T_1,T_2$ is adjacent to $T$ or $T'$, so Proposition~\ref{prop:infinite}(v) applies to give a contradiction.


{Subsubcase 1.2.1.} \emph{There is a vertex $r$ which is not in $Sk$ such that $r$ is connected to $k$.} Then $r$ is connected to exactly one of $s_1,s_3$ by the Claim above. Thus we can assume that $r$ is connected to only $k$ and $s_3$ in $C'$ (if $r$ is connected to $s_4$, then Proposition~\ref{prop:infinite}(iv) applies) and that the triangle $\{k,r,s_3\}$ is oriented (otherwise, since it is adjacent to the non-oriented $C'$, Proposition~\ref{prop:infinite}(v) applies). Note then that $r$ is not connected to $s_2$ or $s_4$ because then it is connected to exactly three vertices in $C$ or $C'$. We will show that $e_{s_2}+e_{s_4}$ is a radical vector. Suppose that this not true, i.e. there is a vertex $t$ which is not in $Sk$ such that (i) $t$ is connected to exactly one of $s_2,s_4$ or (ii) $t$ is connected to both of $s_2,s_4$ such that the edges $\{t,s_2\}$ and $\{t,s_4\}$ have the same orientations (Lemma~\ref{lem:Srad}). We will show that this contradicts our assumptions. Note that since the triangle $\{k,r,s_3\}$ is oriented, any edge incident to $r$ is a single edge (otherwise Proposition~\ref{prop:infinite}(ii) applies to give a contradiction).

We consider the subcases for (i). Let us first suppose that (a) $t$ is connected to $s_4$ and not connected to $s_2$. If $t$ is connected to $k$ as well, then $t$ is connected to exactly one of $s_1,s_3$ by the Claim, then it is connected to exactly three vertices in $C'$, which gives a contradiction. Thus we can assume that $t$ is \emph{not} connected to $k$. Then $t$ is connected to exactly one of $s_1,s_3$ (otherwise $t$ is connected to an odd number of vertices in the non-oriented cycle $C'$). 

Under all these assumptions, suppose that (a1) $t$ is connected to $s_3$ (and not connected to $s_1$). Then note that the triangle $\{t,s_3,s_4\}$ is oriented (otherwise, since it is adjacent to $C'$, Proposition~\ref{prop:infinite}(v) applies). If $t$ is connected to $r$ then the cycle $C''=\{s_4,t,r,k,s_1\}$ is non-oriented (where $s_1$ is a sink) and $s_2$ is connected to exactly one vertex (which is $s_1$) in $C''$, contradiction by Proposition~\ref{prop:infinite}(iv). If $t$ is not connected to $r$, then the subquiver $\{r,k,s_3,s_4,t,s_2\}$ is mutation-equivalent to $E_6$, which is a contradiction. Suppose now that (a2) $t$ is connected to $s_1$ (and not connected to $s_3$). Similarly the triangle $\{t,s_1,s_4\}$ is oriented. If $t$ is connected to $r$ then the cycle $C'''=\{t,r,k,s_1\}$ is non-oriented (where k is a source) and $s_2$ is connected to exactly one vertex (which is $s_1$) in $C'''$, contradiction. If $t$ is not connected to $r$, then the subquiver $\{r,k,s_1,s_4,t,s_2\}$ is mutation-equivalent to $E_6$, contradiction. 

Let us now suppose that (b) $t$ is connected to $s_2$ and not connected to $s_4$.
If $t$ is connected to $k$ as well, then t is connected to exactly one of $s_1,s_3$ by the Claim, so it is connected to exactly three vertices in $C$ contradiction (Proposition~\ref{prop:infinite}(iv)). Thus assume that $t$ is not connected to $k$. (b1) If $t$ is connected to $r$, then we have the following: if $t$ is not connected to $s_1$, then the cycle $C''''=\{k,s_1,s_2,t,r\}$ is a non-oriented cycle (where $k$ is a source) and $s_4$ is connected to exactly one vertex (which is $s_1$) in $C''''$, contradiction; if $t$ is connected to $s_1$, then the cycle $C'''''=\{k,s_1,t,r\}$ is a non-oriented cycle (where $k$ is a source) and $s_4$ is connected to exactly one vertex (which is $s_1$) in $C'''''$, contradiction. (b2) If $t$ is not connected to $r$, then the subquiver $\{t,s_2,s_3,r,k,s_4\}$ is mutation-equivalent to $E_6$, contradiction.

We consider the subcases for (ii), so suppose that $t$ is connected to both $s_2,s_4$ such that the edges $\{t,s_2\}$ and $\{t,s_4\}$ have the same orientations. If $t$ is connected to $k$ as well, then $t$ is connected to exactly one of $s_1,s_3$ by the Claim, so it is connected to exactly three vertices in $C$ contradiction (Proposition~\ref{prop:infinite}(iv)). Thus assume that $t$ is not connected to $k$. Then again $t$ is connected to exactly one of $s_1,s_3$ (otherwise $t$ is connected to an odd number of vertices in the non-oriented cycle $C'$), however then $t$ is connected to exactly three vertices in $S$, so again Proposition~\ref{prop:infinite}(iv) applies to give a contradiction.


{Subsubcase 1.2.2.} \emph{No vertex $r$ which is not in $Sk$ is connected to $k$.}
First suppose that there is a vertex $r$ which is not in $Sk$ such that $r$ is connected to $s_4$. Then (**) $r$ is connected to exactly one of $s_1,s_3$, say connected to $s_1$ (otherwise $r$ is connected to exactly an odd number of vertices in $C'$ because $r$ is not connected to $k$, so Proposition~\ref{prop:infinite}(iv) applies). Then $r$ is not connected to $s_2$ (because then it is connected to exactly three vertices in $S$). This implies, in particular, that the vectors $e_{s_1}+e_{s_3}$ and $e_{s_2}+e_{s_4}$ are both not radical. We claim that $e_{s_1}+e_{s_2}+e_{s_3}+e_{s_4}$ is radical. Suppose that this is not true. Then there is a vertex $t\ne k,r$ connected to $S$ such that, for the edges connecting $t$ to $S$, the number of the ones going away from $t$ is different from the ones going in (Lemma~\ref{lem:Srad}). Then either (i) $t$ is connected to exactly one vertex in $S$ or (ii) $t$ is connected exactly two vertices $s_i,s_j$ and the orientations of the corresponding edges are the same (so $t$ is a source or sink in $St$) (here note that, by (**), the vertex $t$ is not connected to all vertices in $S$, so it is connected to at most two vertices in $S$ by Proposition~\ref{prop:infinite}(iv)).
In the latter case (ii), the vertices $s_i$ and $s_j$ are not connected (because otherwise the triangle $T=\{t,s_i,s_j\}$ is non-oriented and any of the remaining vertices of $S$ is connected to exactly one vertex in $T$), so this case is the same as {Subcase 1.1} above replacing $k$ by $t$, 
which implies that $e_{s_1}+e_{s_3}$ or $e_{s_2}+e_{s_4}$ is radical, however this is not true in this case, so contradiction. 
Thus here we only need to consider the case (i), where $t$ is connected to exactly one vertex, say $v$, in $S$. Then $v=s_2$ because otherwise $t$ is connected to exactly one vertex in the non-oriented cycle $C'$ (note that $t$ is not connected to $k$ by the definition of this case). If $t$ is connected to $r$, then the subquiver $\{t,r,s_1,s_2\}$ is a non-oriented cycle that contains $t$ and $r$
and $k$ is connected to exactly one vertex there (note that the triangle $\{r,s_1,s_4\}$ is oriented because otherwise, since it is adjacent to the non-oriented $C'$, Proposition~\ref{prop:infinite}(v) applies), contradiction. If $t$ is not connected to $r$, then the subquiver $Str$ is mutation-equivalent to $E_6$, which also contradicts an assumption. 


Now suppose that no vertex $r$ which is not in $Sk$ is connected to $s_4$. 
Note that if such a vertex $r$ is connected to $s_1$ or $s_3$, then it is connected to the other one as well (otherwise $t$ is connected to exactly one vertex in $C'$) and it is not connected to $s_2$ (otherwise it is connected to exactly three vertices in $S$ and Proposition~\ref{prop:infinite}(iv) applies). Also if $r$ is connected to both $s_1$ and $s_3$ with the same orientations, then the cycles $\{r,s_1,s_4,s_3\}$ and  $\{r,s_1,s_3,k\}$ are non-oriented and their union has no oriented cycles, so contradiction by Proposition~\ref{prop:infinite}(v). 
Thus if a vertex is connected $s_1$ or $s_3$, it is connected to both of them with opposite orientations. This implies that $e_{s_1}+e_{s_3}$ is radical.



{Case 2.} \emph{$k$ is connected to all four vertices in $S$.} Then $Sk$ has four triangles that contain $k$; two of them are oriented and the other two are non-oriented and the non-oriented ones are not adjacent (Proposition~\ref{prop:infinite}(v)). We will show that the vector $e_{s_1}+e_{s_2}+e_{s_3}+e_{s_4}$ is radical. This trivially follows from Lemma~\ref{lem:Srad} if there is no other vertex which is connected to $S$. Let us now denote by $r$ an arbitrary vertex which is connected to $S$. If $r$ is connected to $k$ by a double edge, then it is connected to an oriented triangle, so Proposition~\ref{prop:infinite}(ii) applies necessarily; if $r$ is connected to $k$ by a single edge, then it is connected to both of $s_1,s_3$ or both $s_2,s_4$, but not all four of them (otherwise $k$ is connected to exactly three vertices in one of the non-oriented triangles of $Sk$ containing $k$, then Proposition~\ref{prop:infinite}(iv) applies) so $Skr$ is mutation-equivalent to $X_6$, contradiction. Let us now assume that $r$ is not connected to $k$ and the subquiver $Skr$ is simply-laced. If $r$ is connected to a vertex in one of the non-oriented triangles, then it is connected to the other vertex as well with the opposite orientation by Proposition~\ref{prop:infinite}(iv,v) (note that every vertex of $S$ lies in exactly one of the non-oriented triangles). Then for the edges connecting $r$ to $S$, the number of the ones leaving $r$ is the same as the ones entering $r$. Thus $e_{s_1}+e_{s_2}+e_{s_3}+e_{s_4}$ is radical by Lemma~\ref{lem:Srad}.

We have completed the proof that the basic quiver which is oriented square contains the support of a non-zero radical vector $u$ as in Theorem~\ref{th:basic radical}. To show this for the other basic subquivers with $4$ vertices, let us first note that they are mutation-equivalent to the oriented square. Furthermore, mutation of a basic quiver with four vertices is also a basic quiver. Let us now suppose that $S$ is a basic subquiver in $Q$ and $u$ is a radical vector as in the theorem. Let $k$ be a vertex in $S$ and let $u'$ be the vector that represents $u$ in the basis which corresponds to ${B}^{Q'}=\mu_k({B}^{Q})$ (Section~\ref{sec:pre}). Then $supp_{Q'}(u')$ lies in the basic subquiver $S'=\mu_k(S)\subset \mu_k(Q)=Q'$. Also $Q'$ does not contain any subquiver which is mutation-equivalent to $E_6$ or $X_6$ (by \cite[Corollary~5.13]{FSTu} or by Theorem~\ref{th:E6 inv}, which we will prove without using the current theorem). Thus we can conclude that any basic subquiver $S$ with four vertices contains the support of a non-zero radical vector as in Theorem~\ref{th:basic radical}. For a basic quiver $S$ with $m>4$ vertices, applying the mutation at a vertex $k$ of $S$ gives a basic subquiver $S'$ with $m-1$ vertices, then the existence of the non-zero radical vector $u$ follows by induction and the base change formula as we discussed.


Conversely, to show the if part of the theorem, suppose that $Q$ is a finite mutation type quiver which does not come from the triangulation of a surface. Then $Q$ contains a subquiver which is mutation-equivalent to $E_6$ or $X_6$ \cite[Corollary~5.13]{FSTu}. In fact, it is enough to show it for $E_6$ and $X_6$ (see also Theorem~\ref{th:E6 inv}). For $E_6$ it follows from Proposition~\ref{prop:E6 char}, for $X_6$ it follows from a direct check using Lemma~\ref{lem:Srad}.


\subsection{Proof of Theorem~\ref{th:V00}}

\label{sec:V00}

The first part of the theorem, where $Q$ contains a double edge or a non-oriented cycle, follows from Proposition~\ref{prop:infinite}(ii,iv) and Lemma~\ref{lem:Srad}. Let us now prove the second part. For this, let us assume that $Q$ is a quiver which comes from the triangulation of a surface or $Q$ is mutation equivalent to one of $X_6,X_7$. 
Then $Q$ does not contain any subquiver which is mutation-equivalent to $E_6$ \cite[Corollary~5.13]{FSTu}. We first show (ii) for convenience. Let $Q'=\mu_k(Q)$. 
Suppose that $u$ be a basic radical vector for $\bar{B}^Q$ and let $u'$ be the vector that represents $u$ in the basis that corresponds to $\bar{B}^{Q'}$ (see Section~\ref{sec:pre} for the base change formula corresponding to $\mu_k$). We will show the following: 

Claim: $u'$ is in the span of basic radical vectors for $\bar{B}^{Q'}$. 

Proof: If $k$ is in $supp_{Q}(u)$, then it follows from a direct check that $u'$ is a basic radical vector for  $\bar{B}^{Q'}$. We consider the case where $k$ is not in $supp_{Q}(u)$ (but connected to it). It is easy to check the claim if $supp_{Q}(u)$ has exactly two vertices. Thus we assume that $C=supp_{Q}(u)$ is a cycle. We denote the subquiver $C\cup \{k\}$ by $Ck$ for convenience.
Note that if $u$ lies in the span of basic radical vectors whose support contain $k$, then we are done. 
First we consider the subcase that $k$ is connected to a vertex, say $c$, in $C$ by a double edge. 
Let $c',c''$ be the vertices which are adjacent to $c$ in $C$. Then, by Proposition~\ref{prop:infinite}(ii), the vertex $k$ is connected to both of $c'$ and $c''$ and it is not connected to any other vertex in $C$. Then the subquiver $C'$ obtained from $Ck$ by removing $c$ is a non-oriented cycle  ($k$ is a source or sink there). Thus $u=y-x$ where $y$ the vector such that $supp_{Q}(y)=C'$ and $x$ is the vector such that $supp_{Q}(x)=\{k,c\}$ (recall that we work modulo $2\ZZ$). the vectors $x,y$ are basic radical vectors for $\bar{B}^{Q}$ by the first part of the theorem. Since the claim is true for $x,y$ (their supports contain $k$), it is also true for $u$. Similar arguments, in view of Theorem~\ref{th:basic radical}, also show the Claim if $Ck$ has a subquiver $S$ which is a non-oriented cycle or a basic subquiver such that $S$ contains $k$. Then the remaining subcase is where the vertex $k$ is connected to exactly two vertices $c_1,c_2$ in $C$ and that the triangle $\{k,c_1,c_2\}$ is oriented. Then $u'$ is a basic radical vector for  $\bar{B}^{Q'}$ with $supp_{Q'}(u)=\mu_k(Ck)$, which is a cycle in $Q'=\mu_k(Q)$.

Thus, by the Claim, we have $dim(\bar{V}^{Q}_{00})\leq dim(\bar{V}^{Q'}_{00})$. Since $\mu_k$ is involutive, it is also true that $dim(\bar{V}^{Q'}_{00})\leq dim(\bar{V}^{Q}_{00})$, so $dim(\bar{V}^{Q}_{00})=dim(\bar{V}^{Q'}_{00})$.

To prove part (i), let us first note that the statement is true if $Q$ is (mutation-equivalent to) the Dynkin quiver $A_n$. For arbitrary $Q$, we will reduce the claim to the $A_n$ case, which will prove the statement. Let us now assume that $u,v$ are two non-zero radical vectors which are not in $\bar{V}^{Q}_{00}$. 
Since we could replace $u$ by $u-w$, we can assume without loss of generality that: 

(****) the union of $supp_Q(u)$ and $supp_Q(v)$ (in particular each of them) does not contain the support of a basic radical vector $w$ of $\bar{B}^Q$.

Then, by Theorem~\ref{th:basic radical}, any connected component of $supp_Q(u)$ (or $supp_Q(v)$) is a single vertex or a simply-laced oriented triangle. In view of part (ii), applying some mutations if necessary, we can assume that each connected component of $supp_Q(u)$ and $supp_Q(v)$ is a single vertex 
(if mutations are applied the conditions of the statement will also be satisfied for the resulting quiver because of Theorem~\ref{th:E6 inv}). Let us now note that (****) implies the following: a minimal connected subquiver $M$ that contains $supp_Q(u)$ and $supp_Q(v)$ does not contain any basic subquiver or a non-oriented cycle or a double edge. This is because if $M$ contains such a subquiver $S$ then, by Theorem~\ref{th:basic radical} and part (ii) of the current theorem, there is a basic radical vector $w$ such that $supp_Q(w)$ lies in $S$; furthermore, by the assumption (****), there is a vertex in $supp_Q(w)$ which is not contained in any of $supp_Q(u)$ or $supp_Q(v)$; removing this vertex gives a connected subquiver (because if a vertex is connected to $supp_Q(w)$ it is connected to at least two vertices there). This  contradicts minimality of $M$. Thus $M$ is mutation-equivalent to $A_n$. Since $u,v$ belong to $V^{M}_{00}$, whose dimension is at most $1$, they are linearly dependent, so $\bar{V}^{Q}_{0}/\bar{V}^{Q}_{00}$ has dimension at most $1$. This completes the proof of the theorem.





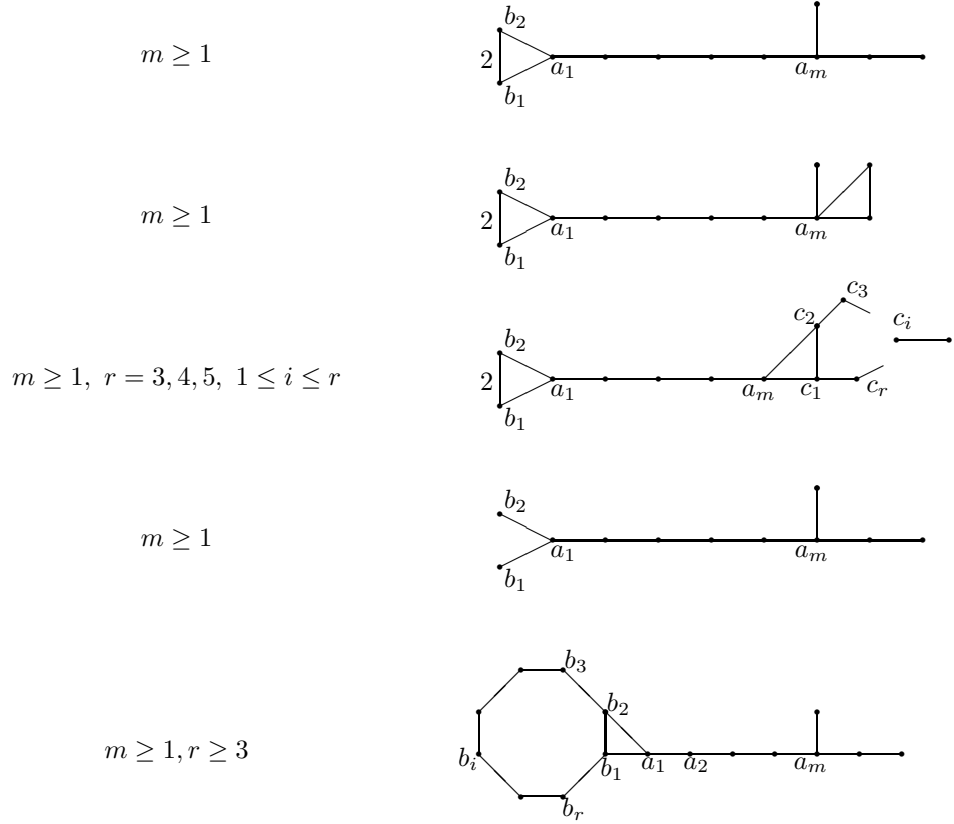
\begin{figure}
\[ 
\begin{array}{ccl} 
{m\geq 1}
&& 
\setlength{\unitlength}{1.0pt}
\begin{picture}(180,43)(-20,-2)
\put(40,0){\circle*{2.0}}
\put(20,10){\circle*{2.0}}  \put(26,15){\makebox(0,0){$b_2$}}
\put(20,-10){\circle*{2.0}} \put(26,-15){\makebox(0,0){$b_1$}}
\put(60,0){\circle*{2.0}}
\put(80,0){\circle*{2.0}}
\put(100,0){\circle*{2.0}}
\put(120,0){\circle*{2.0}}
\put(140,0){\circle*{2.0}}

\put(160,0){\circle*{2.0}}
\put(140,20){\circle*{2.0}}
\put(140,20){\circle*{2.0}}

\put(180,0){\circle*{2.0}}



\put(40,0){\line(1,0){140}}
\put(40,0){\line(-2,-1){20}}
\put(40,0){\line(-2,1){20}}
\put(140,0){\line(0,1){20}}


\put(20,-10){\line(0,1){20}}
\put(15,-1){\makebox(0,0){$2$}}

\put(138,-5){\makebox(0,0){$a_m$}}
\put(44,-5){\makebox(0,0){$a_1$}}

\end{picture}
\\[.2in] 

{m\geq 1}
&& 
\setlength{\unitlength}{1.0pt}
\begin{picture}(180,43)(-20,-2)
\put(40,0){\circle*{2.0}}
\put(20,10){\circle*{2.0}}  \put(26,15){\makebox(0,0){$b_2$}}
\put(20,-10){\circle*{2.0}} \put(26,-15){\makebox(0,0){$b_1$}}
\put(60,0){\circle*{2.0}}
\put(80,0){\circle*{2.0}}
\put(100,0){\circle*{2.0}}
\put(120,0){\circle*{2.0}}
\put(140,0){\circle*{2.0}}

\put(160,0){\circle*{2.0}}
\put(140,20){\circle*{2.0}}
\put(140,20){\circle*{2.0}}

\put(160,20){\circle*{2.0}}



\put(40,0){\line(1,0){120}}
\put(40,0){\line(-2,-1){20}}
\put(40,0){\line(-2,1){20}}
\put(140,0){\line(0,1){20}}

\put(160,0){\line(0,1){20}}
\put(160,20){\line(-1,-1){20}}

\put(20,-10){\line(0,1){20}}
\put(15,-1){\makebox(0,0){$2$}}

\put(138,-5){\makebox(0,0){$a_m$}}
\put(44,-5){\makebox(0,0){$a_1$}}

\end{picture}
\\[.2in] 
{m\geq 1,~r=3,4,5}, ~1\leq i \leq r
&& 
\setlength{\unitlength}{1.0pt}
\begin{picture}(180,43)(-20,-2)
\put(40,0){\circle*{2.0}}
\put(20,10){\circle*{2.0}}  \put(26,15){\makebox(0,0){$b_2$}}
\put(20,-10){\circle*{2.0}} \put(26,-15){\makebox(0,0){$b_1$}}
\put(60,0){\circle*{2.0}}
\put(80,0){\circle*{2.0}}
\put(100,0){\circle*{2.0}}
\put(120,0){\circle*{2.0}}
\put(140,0){\circle*{2.0}}

\put(140,20){\circle*{2.0}}
\put(140,20){\circle*{2.0}}

\put(155,0){\circle*{2.0}}
\put(150,30){\circle*{2.0}}
\put(190,15){\circle*{2.0}}
\put(170,15){\circle*{2.0}}

\put(175,21){\makebox(-5,2){$c_i$}}
\put(138,23){\makebox(-5,2){$c_2$}}
\put(138,-3){\makebox(0,-3){$c_1$}}

\put(158,33){\makebox(-5,2){$c_3$}}
\put(163,-3){\makebox(0,-1){$c_r$}}

\put(40,0){\line(1,0){115}}
\put(40,0){\line(-2,-1){20}}
\put(40,0){\line(-2,1){20}}
\put(140,0){\line(0,1){20}}

\put(120,0){\line(1,1){20}}
\put(140,20){\line(1,1){10}}
\put(150,30){\line(2,-1){10}}
\put(155,0){\line(2,1){10}}
\put(170,15){\line(1,0){20}}

\put(20,-10){\line(0,1){20}}
\put(15,-1){\makebox(0,0){$2$}}

\put(118,-5){\makebox(0,0){$a_m$}}
\put(44,-5){\makebox(0,0){$a_1$}}

\end{picture}
\\[.2in] 

{m\geq 1}
&& 
\setlength{\unitlength}{1.0pt}
\begin{picture}(180,43)(-20,-2)
\put(40,0){\circle*{2.0}}
\put(20,10){\circle*{2.0}}  \put(26,15){\makebox(0,0){$b_2$}}
\put(20,-10){\circle*{2.0}} \put(26,-15){\makebox(0,0){$b_1$}}
\put(60,0){\circle*{2.0}}
\put(80,0){\circle*{2.0}}
\put(100,0){\circle*{2.0}}
\put(120,0){\circle*{2.0}}
\put(140,0){\circle*{2.0}}

\put(160,0){\circle*{2.0}}
\put(140,20){\circle*{2.0}}
\put(140,20){\circle*{2.0}}

\put(180,0){\circle*{2.0}}



\put(40,0){\line(1,0){140}}
\put(40,0){\line(-2,-1){20}}
\put(40,0){\line(-2,1){20}}
\put(140,0){\line(0,1){20}}



\put(138,-5){\makebox(0,0){$a_m$}}
\put(44,-5){\makebox(0,0){$a_1$}}

\end{picture}
\\[.4in] 

{m\geq 1, r\geq 3}
&& 
\setlength{\unitlength}{0.8pt} 
\begin{picture}(140,60)(-20,-2) 
\put(100,0){\circle*{2.0}}

\put(80,20){\circle*{2.0}}
\put(80,0){\circle*{2.0}} 
\put(120,0){\circle*{2.0}}
\put(140,0){\circle*{2.0}}
\put(160,0){\circle*{2.0}}
\put(180,0){\circle*{2.0}}
\put(80,20){\circle*{2.0}}
\put(60,40){\circle*{2.0}}
\put(40,40){\circle*{2.0}}
\put(20,20){\circle*{2.0}}
\put(40,-20){\circle*{2.0}}
\put(20,0){\circle*{2.0}}
\put(60,-20){\circle*{2.0}}

\put(80,0){\line(1,0){140}}
\put(80,0){\line(-1,-1){20}}
\put(80,0){\line(0,1){20}}
\put(100,0){\line(-1,1){20}}
\put(80,20){\line(-1,1){20}}
\put(60,40){\line(-1,0){20}}
\put(40,40){\line(-1,-1){20}}
\put(20,20){\line(0,-1){20}}
\put(20,0){\line(1,-1){20}}
\put(40,-20){\line(1,0){20}}

\put(180,0){\line(0,1){20}}
\put(180,20){\circle*{2.0}}
\put(200,0){\circle*{2.0}}
\put(220,0){\circle*{2.0}}

\put(83,-7){\makebox(0,0){$b_1$}}
\put(66,44){\makebox(0,0){$b_3$}}
\put(15,-2){\makebox(0,0){$b_i$}}
\put(65,-26){\makebox(0,0){$b_r$}}
\put(86,24){\makebox(0,0){$b_2$}}

\put(103,-5){\makebox(0,0){$a_1$}}
\put(123,-5){\makebox(0,0){$a_2$}}
\put(177,-5){\makebox(0,0){$a_m$}}

\end{picture} 
\\[.3in] 
\end{array} 
\] 
\caption{Some quivers of infinite mutation type; each edge can be taken to be arbitrarily oriented (see also Proposition~\ref{prop:infinite}).} 
\label{fig:critical} 
\end{figure} 


\subsection{Proof of Theorem~\ref{th:E6 inv}}

\label{sec:E6 inv}

We prove the theorem in three lemmas:


\begin{lemma}
\label{lem:X6 inv}
Let $Q$ be a quiver of finite mutation type and let $k$ be a vertex in $Q$. Suppose that $Q$ has a subquiver which is mutation-equivalent to $X_6$. Then $\mu_k(Q)$ also contains a subquiver which is mutation-equivalent to $X_6$. 
\end{lemma}

\begin{proof}
Let us denote by $X$ the subquiver which is mutation-equivalent to $X_6$. The lemma is obvious if $k$ is in $X$, so we can assume that $k$ is not in $X$. We then denote the subquiver $\{X,k\}$ by $Xk$ for convenience. 
For quivers which are mutation-equivalent to $X_7$ the lemma follows from a direct check (these quivers are given in \cite{DO}). Then, to complete the proof of the lemma for a general $Q$, it is enough to show that $Xk$ is mutation-equivalent to $X_7$. For this purpose, applying mutations if necessary, we can assume that $X=X_6$. Let us then denote the double edges of $X$ by $\{i_1,i_2\}$ and $\{j_1,j_2\}$, the center vertex by $c$ and the remaining vertex by $d$. By Proposition~\ref{prop:infinite}(i,ii), we can assume without loss of generality that $k$ is connected to $\{i_1,i_2\}$ such that the triangle $\{i_1,i_2,k\}$ is oriented. If $k$ is not connected to any other vertex in $X_6$, then the subquiver obtained from $Xk$ by removing $i_1$ (or $i_2$) is of infinite mutation type (it belongs to Figure~\ref{fig:critical}), contradicting that $Q$ is of finite mutation type. If $k$ is connected to a vertex which is different from $i_1,i_2$, then it follows from a direct check that there is a non-oriented cycle that contains $k$ such that Proposition~\ref{prop:infinite}(iv,v) applies. Thus $Xk$ is mutation-equivalent to $X_7$. This completes the proof of the lemma.
\end{proof}

\begin{lemma}
\label{lem:E6 inv}
Suppose that $Q$ is a quiver of finite mutation type and let $k$ be a vertex in $Q$. Suppose also that $Q$ has a subquiver which is mutation-equivalent to $E_6$. Then $Q$ does not contain any subquiver which is mutation equivalent to $X_6$. Furthermore, the quiver $\mu_k(Q)=Q'$ also contains a subquiver which is mutation-equivalent to $E_6$. 
\end{lemma}

\begin{proof} 
For the first part of the lemma, let us assume to the contrary that $Q$ contains a subquiver which is mutation-equivalent to $X_6$. Then, by Lemma~\ref{lem:X6 inv}, the quiver $Q$ is mutation-equivalent to $X_6$ or $X_7$. However, it follows from a direct check on the quivers which are mutation equivalent to $X_6$ or $X_7$ as given in \cite{DO}, that $Q$ does not contain any subquiver which is mutation-equivalent to $E_6$, contradicting our assumption.

For the second part of the lemma, we first note that, by \cite[Theorem~5.11]{FSTu}, the quiver $Q$, so $Q'$, does not come from the triangulation of a surface. Therefore $Q'$ also contains a subquiver which is mutation eqivalent to $E_6$ or $X_6$ (which are minimal quivers that do not come from the triangulation of a surface). By the first part of the lemma, the quiver $Q'$ does not contain any subquiver which is mutation-equivalent to $X_6$, thus $Q'$ contains a subquiver which is mutation-equivalent to $E_6$. This completes the proof of the lemma.
\end{proof}

\begin{lemma}
\label{lem:E6-classify}
Let $Q$ be a connected quiver of finite mutation type. Suppose also that $Q$ has a subquiver which is mutation-equivalent to $E_6$. Then $Q$ is mutation-equivalent to a quiver which is one of the (exceptional) types $E_6$, $E_7$, $E_8$, $E_6^{(1)}$, $E_7^{(1)}$, $E_8^{(1)}$, $E_7^{(1,1)}$, $E_8^{(1,1)}$ \cite[Figure~6.1]{FSTu}.
\end{lemma}

\begin{proof} 
If $Q$ is of finite type, then the lemma follows from the classification of finite type quivers under the mutation operation \cite{CAII}. Thus we can assume that $Q$ is not of finite type. Then it is mutation-equivalent to a quiver $Q'$ which has a double edge $e=\{u_1,u_2\}$. Note that the quiver $Q'$ also contains a subquiver which is mutation-equivalent to $E_6$ by Lemma~\ref{lem:E6 inv}. Below we consider cases depending on the number of vertices connected to $e$. In our analysis, if we do not specify an orientation on a subquiver, it is assumed to be oriented as required by 
Proposition~\ref{prop:infinite} to be of finite mutation type.

Case 1. \emph{There is exactly one vertex, say $v_1$, connected to $e$}. 
In this case it can be checked easily, using Proposition~\ref{prop:E6 char}(iii), that $Q'$ contains a subquiver $M'$ as in Figure~\ref{fig:critical}. (More explicitly $M'$ is a minimal connected subquiver which contains $e$ and a subquiver $E'$ which is mutation-equivalent to $E_6$). Since the quivers in Figure~\ref{fig:critical} are of infinite mutation type, the quiver $Q'$ is also of infinite mutation type, contradicting our assumption.



Case 2. \emph{There are exactly two vertices, say $v_1,v_2$, connected to $e$}. By Proposition~\ref{prop:infinite}(ii,iii), we can assume that the subquiver $\{e,v_i\}$ is an oriented triangle for $i=1,2$. Let us also note that if there is a subquiver $E'$ which is mutation-equivalent to $E_6$ such that $E'$ contains at most one of $v_1$ and $v_2$, then $Q'$ contains a subquiver as in Case 1 (it is a minimal connected subquiver that contains $E'$ together with the edge $e$ and one of the vertices $v_1$ or $v_2$), which we already considered. Thus here we only need to consider the case when any subquiver $E'$ which is mutation-equivalent to $E_6$ contains both $v_1$ and $v_2$. 
Let us note that, by Proposition~\ref{prop:E6 char}(iii), there are vertices $v'_1$ and $v'_2$ in $E'$ which are connected to $v_1$ and $v_2$ respectively. Let us denote $e_1=\{v_1,v'_1\}$, $e_2=\{v_2,v'_2\}$ 

We first consider the subcase where $E'$ contains a vertex which is adjacent to $e$ (i.e. $E'$ contains $u_1$ or $u_2$).
Then any path connecting any of $v_1,v'_1$ to $v_2$ or $v'_2$ in $E'$ contains $u_1$ or $u_2$ because otherwise there is a non-oriented cycle in $E'$ (note that, since the subquiver $\{e,v_i\}$ is an oriented triangle, the edges that connect $v_1$ and $v_2$ to $u_j, j=1,2$ have the same orientation). In particular, we have $v'_1\ne v'_2$. Furthermore, the remaining vertex of $E'$ is connected to exactly one of  $e_1$ and $e_2$. However, then $E'$ is either the Dynkin tree $D_6$ or it is mutation-equivalent to $A_6$,
contradiction.

Let us now consider the subcase where $E'$ does not contain any vertex which is in $e$ (i.e. $E'$ does not contain any of $u_1,u_2$). By Proposition~\ref{prop:E6 char}, $E'$ contains a basic subquiver $S$. Since $S$ has at least four vertices, it contains at least two vertices from the set $\{v_1,v'_1,v_2,v'_2\}$. Let us first assume that $S$ does not contain any vertex from one of the edges $e_1,e_2$, say it does not contain any vertex from $e_2$. Then $S$ has exacly four vertices (because $E'$ has six vertices and it already contains $e_1$ and $e_2$). Also the vertices in $e_2$ are connected to $S$ as required Proposition~\ref{prop:E6 char}(ii,iii). Then, it follows from an easy check that $Q'$ contains a non-oriented cycle $C$ such that there is a vertex $r$ (in $E'$) which is connected (by singles edges) to an odd number of vertices in $C$, contradiction by Proposition~\ref{prop:infinite}(iv).

To complete the treatment of this case, we now assume that $S$ contains a vertex from each of $e_1$ and $e_2$ and consider the subcases:

Subcase 2.1. \emph{$S$ is the Dynkin tree $D_4$}. Let $S=\{a,b,c,d\}$ such that $a$ is in $e_1$; $b$ is in $e_2$; $c$ is connected to both $a,b$ and $d$ is connected to only $c$ in $S$ (so $c$ is the "center" of $S$). Let $a',b'$ be the remaining vertices in $e_1$ and $e_2$ respectively (so $\{a,a'\}=\{v_1,v'_1\}$ and $\{b,b'\}=\{v_2,v'_2\}$). 

Let us first assume that no vertex in $e_1$ is connected to any vertex in $e_2$. Then the vertex $d$ is not connected to $a'$ nor to $b'$, by Proposition~\ref{prop:E6 char}(iii) (applied to the pairs of vertices $d,a$ and $d,b$). Let us note that, since the triangles $\{e,v_1\}$ and $\{e,v_2\}$ are oriented, there are two non-oriented cycles in $Q'$ which contain $v_1,v_2$ and the vertex $c$. The vertex $d$ is connected to exactly one vertex (which is $c$) in these cycles, contradiction by Proposition~\ref{prop:infinite}(iv). If there is a vertex in $e_1$ which is connected to a vertex in $e_2$, the subcase follows by similar arguments.

Subcase 2.2. \emph{$S$ is formed by two adjacent triangles}. As in the previous subcase, let $S=\{a,b,c,d\}$ such that $a$ is in $e_1$; $b$ is in $e_2$; $c,d$ are connected to both $a,b$ and to each other (so the triangles of $S$ are $\{a,c,d\}$ and $\{b,c,d\}$). Let $a',b'$ be the remaining vertices in $e_1,e_2$ respectively (so $\{a,a'\}=\{v_1,v'_1\}$ and $\{b,b'\}=\{v_2,v'_2\}$).
If there is a vertex in $e_1$ which is connected to a vertex in $e_2$, then there is a non-oriented cycle in $E'$, which is not the case, so we can assume that no vertex in $e_1$ is connected to any vertex in $e_2$. Then, by Proposition~\ref{prop:E6 char}(ii), the vertices $c,d$ are not connected to $a'$ or $b'$. Since the triangles $\{e,v_1\}$ and $\{e,v_2\}$ are oriented, there are two non-oriented cycles say $C_1,C_2$ in $Q'$ which contain $v_1,v_2$ together with the vertex $c$ or $d$ respectively. Then, e.g., the vertex $c$ is connected to exactly three vertices in $C_2$, contradiction by Proposition~\ref{prop:infinite}(iv).

Subcase 2.3. \emph{$S$ is a square}. Let $S=\{a,b,c,d\}$ such that $a$ is in $e_1$, $b$ is in $e_2$, the vertices $c,d$ are connected to both $a,b$ and not connected to each other (so $\{a,a'\}=\{v_1,v'_1\}$ and $\{b,b'\}=\{v_2,v'_2\}$). 
As in the previous subcase, if there is a vertex in $e_1$ which is connected to a vertex in $e_2$, then there is a non-oriented cycle in $E'$, which is not the case, so we can assume that no vertex in $e_1$ is connected to any vertex in $e_2$.
By Proposition~\ref{prop:E6 char}(iii), one of the vertices $a',b'$, say $a'$, is connected to exactly one of $c$ or $d$, say connected to $c$. As in the previous subcases, there are two non-oriented cycles say $C_1,C_2$ in $Q'$ which contain $v_1,v_2$ together with the vertex $c$ or $d$. Then $a'$ or $c$ is connected to exactly an odd number of vertices in one of these non-oriented cycles, contradiction by Proposition~\ref{prop:infinite}(iv).

In the current set-up of this case, it is not possible that $S$ is a cycle with five vertices, so we have completed our analysis for this case.


Case 3. \emph{There are exactly three vertices, say $v_1,v_2,v_3$, connected to $e$}. Let us first note that if there is a subquiver $E'$ which is mutation-equivalent to $E_6$ such that $E'$ contains at most two of $v_1,v_2$, then we are in Case 2. Thus here we only need to consider the case when any  subquiver $E'$ which is mutation-equivalent to $E_6$ contains all $v_1,v_2,v_3$. For any $v_i$, $i=1,2,3$ we denote by $Pv_i$ the subquiver on the vertices which are connected to $v_i$ by a path that does not contain any vertex adjacent to $e$ ($v_i$ is included in $Pv_i$). 
We first show that for any $v_i\ne v_j$ connected to $e$, the subquivers $Pv_i$ and $Pv_j$ are disjoint. Suppose this is not true and assume without loss of generality that $i=1$, $j=2$. Then there is a path $P'=\{v_1=w_1,w_2,...,w_r=v_2\}$, $r\geq 2$, that connects $v_1$ and $v_2$ such that $P'$ does not contain any of $u_1,u_2$ (which are the vertices adjacent to $e$). We can assume without loss of generality that $P'$ a shortest path connecting two vertices which are connected to $e$, implying that $v_3$ is not connected to any vertex in $P'$ except possibly to $v_1$ or $v_2$. If $v_3$ is not connected any of $v_1,v_2$, then the cycle $C=\{P',u_1\}$ is non-oriented (because $u_1$ is connected to $v_1$ and $v_2$ by the same orientation) and $v_3$ is connected to exactly one vertex (which is $u_1$) in $C$, which is a contradiction by Proposition~\ref{prop:infinite}(iv). If $v_3$ is connected any of $v_1,v_2$ then similarly Proposition~\ref{prop:infinite}(iv) or (v) applies to give a contradiction. Thus for the rest of this case, we can assume that the subquivers $Pv_i$ and $Pv_j$ are disjoint.

Let us now note that by Proposition~\ref{prop:E6 char}(iii) at least two of $Pv_i$'s, say $Pv_1$ and $Pv_2$, have at least two vertices. This implies that each $Pv_i$ does not contain any subquiver which is one of the following: a basic subquiver, a non-oriented cycle or a double edge, because otherwise $Q'$ contains a subquiver as in Figure~\ref{fig:critical}. Thus each $Pv_i$ is mutation-equivalent to the Dynkin quiver $A_n$, applying some mutations if necessary we can assume that each $Pv_i$ is of type $A_n$ such that $v_i$ is an end vertex of $Pv_i$ (otherwise $Q'$ also contains a subquiver as in Figure~\ref{fig:critical}). Now we can proceed to establish the lemma:

(i) Suppose that each $Pv_i$ has at at least two vertices. If each of them has exactly two vertices, then $Q'$ is mutation-equivalent to $E_6^{(1,1)}$; if one of them has more than two vertices, then $Q'$ contains a tree which is extended Dynkin (it contains $E_6^{(1)}$ as a proper subquiver), so it is of infinite mutation type \cite{BR}.

Thus for the rest of the proof we can assume that $Pv_3$ has exactly one vertex. Then we have the following subcases:

(ii) Suppose that each of $Pv_1$ and $Pv_2$ has at least three vertices. If both have exactly two vertices, then $Q'$ is mutation-equivalent to $E_7^{(1,1)}$; otherwise $Q'$ contains a tree which is not extended Dynkin (it contains $E_7^{(1)}$ as a proper subquiver), so it is of infinite mutation type, which is a contradiction.

(iii) Suppose now, without loss of generality, that $Pv_2$ has exactly two vertices. If $Pv_1$ has exactly two vertices then $Q'$ is mutation-equivalent to $E_6^{(1)}$; if $Pv_1$ has exactly three vertices, then $Q'$ is mutation-equivalent to $E_7^{(1)}$; if $Pv_1$ has exactly four vertices, then $Q'$ is mutation-equivalent to $E_8^{(1)}$; if $Pv_1$ has exactly five vertices, then $Q'$ is mutation-equivalent to $E_8^{(1,1)}$; if $Pv_1$ has more than five vertices, then $Q'$ contains a tree which is not extended Dynkin (it contains $E_8^{(1)}$ as a proper subquiver),  so it is of infinite mutation type, which is a contradiction. 

Case 4. There are at least four vertices connected to $e$. Let us assume that $v_1,v_2,v_3,v_4$ are connected to $e$. Then the subquiver $S=\{u_1,v_1,v_2,v_3,v_4\}$ is the extended Dynkin tree $D^{(1)}_4$. Since $Q'$ contains a subquiver which is mutation-equivalent to $E_6$, there is a vertex which is connected to $e$ or $S$. Then there is necessarily a tree that contains $S$ as a proper subquiver, so it is of infinite mutation type, which is a contradiction.

This completes the proof of the lemma.
\end{proof} 

Given these lemmas, let us now show how Theorem~\ref{th:E6 inv} follows. Let us first assume that $Q$ is a finite mutation type quiver which contains a subquiver which is mutation equivalent to $E_6$. Then, by Lemma~\ref{lem:E6 inv}, any quiver which is mutation-equivalent to $Q$ contains a subquiver which is mutation-equivalent to $E_6$; furthermore $Q$ is mutation equivalent to a quiver which is one of the (exceptional) types $E_6$, $E_7$, $E_8$, $E_6^{(1)}$, $E_7^{(1)}$, $E_8^{(1)}$, $E_7^{(1,1)}$, $E_8^{(1,1)}$ by Lemma~\ref{lem:E6-classify}. Let us now assume that $Q$ is a finite mutation type quiver which contains a subquiver which is mutation-equivalent to $X_6$. Then by Lemma~\ref{lem:X6 inv} and its proof, any quiver which is mutation-equivalent to $Q$ contains a subquiver which is mutation-equivalent to $X_6$ and $Q$ is in fact mutation-equivalent to the quiver $X_6$ or $X_7$. This completes the proof of the theorem.

\end{document}